\documentclass[12pt,oneside,openany,article]{memoir}
\usepackage{mempatch}
\nouppercaseheads 
\usepackage[amsthm,thmmarks,hyperref]{ntheorem}
\usepackage{xr-hyper}
\usepackage[noTeX]{mmap}
\usepackage{etex}
\usepackage[english]{babel}
\usepackage[leqno]{mathtools}
\usepackage{mathrsfs}
\usepackage{upgreek}
\usepackage[compress,square,comma,numbers]{natbib}
\usepackage{hypernat} 
\usepackage{sfmath}
\usepackage{microtype} 
\usepackage{hyperxmp}
\usepackage{zref}
\usepackage{nameref}
\usepackage[pdftex,bookmarks,pdfnewwindow,plainpages=false,unicode]{hyperref}
\usepackage{bookmark}
\usepackage{attachfile}
\usepackage{enumitem}
\usepackage{tikz}
\usepackage{amssymb} 

\usepackage[color]{showkeys}
\definecolor{refkey}{rgb}{0,0,1}
\definecolor{labelkey}{rgb}{1,0,0}

\hypersetup{
colorlinks=true,
linkcolor=black,
citecolor=black,
urlcolor=blue,
pdfauthor={Victor Ivrii},
pdftitle={Asymptotics of the ground state energy for atoms and molecules in  the self-generated magnetic field},
pdfsubject={Sharp Spectral Asymptotics},
pdfkeywords={Schrodinger operator, self-generated magnetic field},
bookmarksdepth={4}
}

\numberwithin{equation}{chapter}

\theoremstyle{plain}
\newtheorem{theorem}{Theorem}[chapter]

\newtheorem{corollary}[theorem]{Corollary}

\theoremstyle{definition}

\theoremstyle{remark}
\newtheorem{remark}[theorem]{Remark}

\newtheorem{problem}[theorem]{Problem}

\DeclareMathAlphabet{\mathpzc}{OT1}{pzc}{m}{it}

 \newcommand{\cE}{\mathcal{E}}

 \newcommand{\sH}{\mathscr{H}}

 \newcommand{\sL}{\mathscr{L}}

\newcommand{\TF}{{\mathsf{TF}}}

\newcommand{\D}{{\mathsf{D}}}
\newcommand{\E}{{\mathsf{E}}}

\newcommand{\x}{{\mathsf{x}}}

\newcommand{\bC}{{\mathbb{C}}}

\newcommand{\bR}{{\mathbb{R}}}

\newcommand{\fH}{{\mathfrak{H}}}

\def\1{\boldsymbol {|}}
\newcommand{\boldsigma}{{\boldsymbol{\sigma}}}

\newcommand{\blangle}{{\boldsymbol{\langle}}}
\newcommand{\brangle}{{\boldsymbol{\rangle}}}

\newcommand{\Def}{\mathrel{\mathop:}=}

\newcommand{\Spec}{\operatorname{Spec}}

\newcommand{\Tr}{\operatorname{Tr}}

\newenvironment{claim}[1][{\textup{(\theequation)}}]{\refstepcounter{equation}\vglue10pt
\begin{trivlist}
\item[{\hskip\labelsep#1}]}{\vglue10pt\end{trivlist}}

\newenvironment{claim*}[1][{}]{\vglue10pt
\begin{trivlist}
\item[{\hskip\labelsep#1}]}{\vglue10pt\end{trivlist}}

\newcounter{note}

\DeclareTextCommand{\textinfty}{PU}{\9042\036}

\DeclareTextCommand{\textge}{PU}{\9042\145}
\DeclareTextCommand{\textle}{PU}{\9042\144}
\DeclareTextCommand{\texthat}{PD1}{\136}

\begin{document}


\title{Asymptotics of the ground state energy for atoms and molecules in  the self-generated magnetic field}
\author{Victor Ivrii\thanks{Department of Mathematics, University of Toronto, 40 St. George Street, Toronto, ON, M5S 2E4, Canada, ivrii@math.toronto.edu.}}

\maketitle

\chapter{Problem}
\label{sect-1}

This is a last in the series of three papers (following \cite{MQT10, MQT11}) and the theorem~\ref{thm-1-1} and corollary~\ref{cor-1-2} below constitute the final goal of this series. Arguments of this paper are rather standard; all the heavy lifting was done before. Let us consider the following operator (quantum Hamiltonian)

\begin{gather}
\mathsf{H}=\sum_{1\le j\le N}H ^0_{x_j}+\sum_{1\le j<k\le N} |x_j-x_k| ^{-1}
\label{1-1}\\
\shortintertext{in}
\fH= \bigwedge_{1\le n\le N} \sH, \qquad \sH=\sL^2 (\bR^3, \bC^2)\label{1-2}\\
\shortintertext{with}
H ^0=\bigl((i\nabla -A)\cdot \upsigma \bigr) ^2-V(x)
\label{1-3}
\end{gather}
Let us assume that
\begin{claim}\label{1-4}
Operator $\mathsf{H}$ is self-adjoint on $\fH$.
\end{claim}
We will never discuss this assumption. We are interested in the \emph{ground state energy\/} $\E^*_N (A)$ of our system i.e. in the lowest eigenvalue of the operator $\mathsf{H}$ on $\fH$:
\begin{gather}
\E ^*_{N} (0)=\inf \Spec \mathsf{H} \qquad \text{on\ \ } \fH 
\label{1-5}\\
\intertext{as $A=0$ and more generally in} 
\E ^*_N=\inf_A \Bigl( \inf \Spec_{\fH}\mathsf{H} + 
\frac{1}{\alpha} \int |\nabla \times A|^2\,dx\Bigr)
\label{1-6}\\
\shortintertext{where} 
V(x)=\sum_{1\le m\le M} \frac{Z_m}{|x-\x_m|}
\label{1-7}\\
N \approx Z\gg 1, \qquad Z\Def  Z_1+\ldots Z_M, \qquad  Z_1 > 0,\ldots, Z_M>0\label{1-8}\\
\intertext{$M$ is fixed, under assumption}
0<\alpha \le \kappa^*Z^{-1}\label{1-9}
\end{gather}
with sufficiently small constant $\kappa^*>0$.

Our purpose is to prove

\begin{theorem}\label{thm-1-1}
Under assumption \textup{(\ref{1-9})} as $N\ge Z- CZ^{-\frac{2}{3}}$
\begin{gather}
\E^*_N = \cE^\TF_N + \sum_{1\le m\le M} 2Z_m^2 S(\alpha Z_m) + O\bigl(N^{\frac{16}{9}}+ \alpha a^{-3}N ^2 \bigr)\label{1-10}\\
\shortintertext{provided}
a\Def \min _{1\le m <m'\le M} |\x_m-\x_{m'}| \ge N^{-\frac{1}{3}}
\label{1-11}
\end{gather}
where $\cE^\TF_N$ is a \emph{Thomas-Fermi energy\/} (see \textup{\cite{Lieb-1}} or \textup{\cite{ivrii:ground}}) and $S(Z_m)Z_m^2$ are magnetic \emph{Scott correction terms\/} (see \textup{\cite{EFS3}} or \textup{\cite{MQT11}}).
\end{theorem}

Combining with the properties of the Thomas-Fermi energy we arrive to 

\begin{corollary}\label{cor-1-2}
Let us consider $\x_m=\x_m^0$ minimizing full energy 
\begin{gather}
\E^*_N +\sum _{1\le m <m'\le M} Z_mZ_{m'}|\x_m-\x_{m'}|^{-1}.
\label{1-12}\\
\shortintertext{Assume that}
Z_m \asymp N\qquad \forall m=1,\ldots,M.
\label{1-13}
\end{gather}
Then $a\ge N^{-\frac{1}{4}}$ and the remainder estimate in \textup{(\ref{1-10})} is $O\bigl(N^{\frac{16}{9}}\bigr)$.
\end{corollary}

\begin{remark}\label{rem-1-3}
As $\alpha=0$ the remainder estimate (\ref{1-12}) was proven in \cite{ivrii:ground} and the remainder estimate $O\bigl(N^{\frac{5}{3}}(N^{-\delta}+ a^{-\delta})\bigr)$ in \cite{FS} for atoms ($M=1$) and \cite{ivrii:MQT1} for $M\ge 1$; this better asymptotics contains also Dirac and Schwinger correction terms. Unfortunately I was not able to recover such remainder estimate here unless $\alpha$ satisfies stronger assumption than (\ref{1-9}). I still hope to achieve this better estimate without extra assumptions.
\end{remark}

Recall that Thomas-Fermi potential $W^\TF$ and Thomas-Fermi density $\rho^\TF$ satisfy equations
\begin{gather}
\rho^\TF = \frac{1}{3\pi^2} (W^\TF)^{\frac{3}{2}}\label{1-14}\\
\shortintertext{and}
W^\TF =V^0 + \frac{1}{4\pi}|x|^{-1}* \rho^\TF.\label{1-15}
\end{gather}

We prove theorem~\ref{thm-1-1} in sections~\ref{sect-2} ``\nameref{sect-2}'' and~\ref{sect-3} ``\nameref{sect-3}''. Section~\ref{sect-4}  ``\nameref{sect-4}''is devoted to corollary~\ref{cor-1-2} and a brief discussion.

\chapter{Lower estimate}
\label{sect-2}

Consider corresponding to $\mathsf{H}$ quadratic form
\begin{multline}
\blangle \mathsf{H} \Psi, \Psi \brangle = \sum_j (H^0_{x_j}\Psi,\Psi) +  
(\sum_{1\le j<k\le N} |x_j-x_k| ^{-1}\Psi,\Psi)=\\
\sum_j (H_{x_j}\Psi,\Psi) + ((V-W)\Psi,\Psi) +
(\sum_{1\le j<k\le N} |x_j-x_k| ^{-1}\Psi,\Psi)
\label{2-1}
\end{multline}
with
\begin{equation}
H =\bigl((i\nabla -\mathbf {A})\cdot \boldsigma \bigr) ^2-W(x)
\label{2-2}
\end{equation}
where we select $W$ later. By Lieb-Oxford inequality the last term is estimated from below:
\begin{gather}
\blangle \sum_{1\le j<k\le N} |x_j-x_k| ^{-1}\Psi,\Psi\brangle \ge \D(\rho_\Psi,\rho_\Psi) - 
C\int \rho_\Psi^{\frac{4}{3}}\,dx
\label{2-3}\\
\shortintertext{where}
\rho_\Psi (x)= N\int |\Psi(x;x_2,\ldots,x_N)|^2\,dx_2\cdots dx_N
\label{2-4}
\intertext{is a spatial density associated with $\Psi$ and}
\D(\rho,\rho')\Def \frac{1}{2}\iint |x-y|^{-1}\rho(x)\rho'(y)\,dxdy
\label{2-5}
\end{gather}
Therefore
\begin{multline}
\blangle \mathsf{H} \Psi, \Psi \brangle \ge \\
\sum_j (H_{x_j}\Psi,\Psi) - 2((V-W)\Psi,\Psi) +
\D(\rho_\Psi,\rho_\Psi) - 
C\int \rho_\Psi^{\frac{4}{3}}\,dx=\\
\shoveright{\sum_j (H_{x_j}\Psi,\Psi) - 2\D(\rho,\rho_\Psi) +
\D(\rho_\Psi,\rho_\Psi) - 
C\int \rho_\Psi^{\frac{4}{3}}\,dx=\quad\ }\\
\sum_j (H_{x_j}\Psi,\Psi) - \D(\rho,\rho) +
\D(\rho-\rho_\Psi,\rho-\rho_\Psi) - 
C\int \rho_\Psi^{\frac{4}{3}}\,dx
\label{2-6}
\end{multline}
as
\begin{equation}
W-V= |x|^{-1}*\rho.
\label{2-7}
\end{equation}

Note that due to antisymmetry of $\Psi$ 
\begin{equation}
\sum_j (H_{x_j}\Psi,\Psi)\ge  \sum_{1\le j\le N: \lambda_j<0} \lambda_j\ge \Tr^- (H)
\label{2-8}
\end{equation}
where $\lambda_j$ are eigenvalues of $H$.

To estimate the last term in (\ref{2-6}) we reproduce the proof of Lemma \ref{ES3-lm:lo} from \cite{ES3}:

According to magnetic Lieb-Thirring inequality for $U\ge 0$:
\begin{equation}
\sum_{j\le N} \blangle (H^0_{x_j} -U)\Psi,\Psi \brangle \ge
-C\int U^{5/2}\,dx -C\gamma^{-3} U^4\,dx -\gamma \int \mathsf{B}^2dx
\label{2-9}
\end{equation}
$\mathsf{B}=\nabla \times \mathsf{A}$, $\gamma>0$ is arbitrary. Selecting 
$U=\beta \min (\rho_\Psi^{5/3},\gamma \rho_\Psi^{4/3})$ with $\beta>0$ small but independent from $\gamma$ we ensure
$\frac{1}{2}U\rho_\Psi \ge CU^{5/2}+C\gamma^{-3}U^4$ and then
\begin{equation}
\sum_{j\le N}  \blangle (H^0_{x_j})\Psi,\Psi \brangle \ge
\epsilon \int \min (\rho_\Psi^{5/3} , \gamma \rho^{4/3})dx
-\gamma \int \mathsf{B}^2\,dx
\label{2-10}
\end{equation}
which implies 
\begin{multline}
\int \rho_\Psi^{4/3}dx \le 
\gamma^{-1} \int \min (\rho_\Psi^{5/3} , \gamma \rho^{4/3})dx + 
\gamma \int \rho_\Psi dx \le \\
c\gamma^{-1} \sum_{j:\lambda_j<0} \blangle (H^0_{x_j})\Psi,\Psi \brangle +
c \int \mathsf{B}^2 dx + 
c \gamma N
\label{2-11}
\end{multline}
where we use
$\int \rho_\Psi dx =N$. 

\begin{remark}\label{rem-2-1}
As one can prove easily (see also \cite{ES3}) that 
\begin{equation}
\sum_{j\le N} \blangle (H^0_{x_j})\Psi,\Psi \brangle \le CZ^{\frac{4}{3}}N
\label{2-12}
\end{equation}
we conclude that
\begin{equation}
\int \rho_\Psi^{4/3}dx \le CZ^{\frac{2}{3}}N + C_1\int \mathsf{B}^2 dx.
\label{2-13}
\end{equation}
It is sufficient unless we want to recover Dirac-Schwinger terms which unfortunately are too far away for us.
\end{remark}

Therefore skipping the non-negative third term in the right-hand expression of (\ref{2-6}) we conclude that
\begin{multline}
\blangle \mathsf{H} \Psi, \Psi \brangle +
 \frac{1}{\alpha} \int |\nabla \times A|^2\,dx \ge \\
\Tr^-(H) + (\frac{1}{\alpha} -C_1) \int |\nabla \times A|^2\,dx
- \D(\rho,\rho) - CN^{\frac{5}{3}}.
\label{2-14}
\end{multline}
Applying Theorem~\ref{MQT11-thm-5-2} from \cite{MQT11} we conclude that 
\begin{claim}\label{2-15}
the sum of the first and the second terms in the right-hand expression of (\ref{2-14}) is greater than
\begin{equation*}
\frac{2}{15\pi^2} \int W^{\frac{5}{2}}\,dx +\sum_m 2Z_m^2S(\alpha Z_m)- 
C N^{\frac{16}{9}}- C\alpha a^{-3}N^2.
\end{equation*}
\end{claim}
\vspace{-10pt}
To prove this estimate one needs just to rescale $x\mapsto xN^{\frac{1}{3}}$, 
$a\mapsto aN^{\frac{1}{3}}$ and introduce $h=N^{-\frac{1}{3}}$ and $\kappa=\alpha N $. Here one definitely needs the regularity properties like in \cite{MQT11} but we have them as $\rho=\rho^\TF$, $W=W^\TF$. Also one can see easily that ``$-C_1$'' brings correction  not exceeding $C_2\alpha N^2$ as $\alpha N\le 1$. 

Meanwhile  for $\rho=\rho^\TF$, $W=W^\TF$ 
\begin{equation}
\frac{2}{15\pi^2} \int W^{\frac{5}{2}}\,dx -\D(\rho,\rho) =\cE^\TF.
\label{2-16}
\end{equation}
Lower estimate of Theorem~\ref{thm-1-1} has been proven.

\begin{remark}\label{rem-2-2}
$\rho=\rho^\TF$, $W=W^\TF$ delivers the maximum of the right-hand expression of (\ref{2-16}) among $\rho, W$ satisfying (\ref{2-7}).
\end{remark}

\chapter{Upper Estimate}
\label{sect-3}

Upper estimate is easy. Plugging as $\Psi$ the \emph{Slater determinant\/} (see~\cite{ivrii:ground} f.e.)  of $\psi_1,\ldots, \psi_N$ where $\psi_1,\ldots,\psi_N$ are eigenfunctions of $H_{A,W}$  we get 
\begin{multline}
\blangle \mathsf{H} \Psi, \Psi \brangle =\Tr^- (H_{A,W}-\lambda_N) +\lambda_N N + \\
\int (W-V)(x)\rho_\Psi(x)\, dx +
\D (\rho_\Psi, \rho_\Psi ) -\\
\frac{1}{2}N(N-1)
\iint |x_1-x_2|^{-1} |\Psi (x_1,x_2;x_3,\ldots,x_N )|^2\,dx_1\cdots dx_N
\label{3-1}
\end{multline}
where we don't care about last term as we drop it (again because we cannot get sharp enough estimate)  and the first term in the second line is in fact
\begin{equation}
- 2\D(\rho  ,\rho_\Psi );
\label{3-2}
\end{equation}
provided (\ref{2-7}) holds. Thus we get
\begin{multline}
\Tr^- (H_{A,W}-\lambda_N) +\lambda_N N -\D(\rho ,\rho) +
\D (\rho_\Psi -\rho  , \rho_\Psi -\rho) +\\
\frac{1}{\kappa} \int |\partial A|^2\,dx
\label{3-3}
\end{multline}
where we added magnetic energy. Definitely we have several problems here: $\lambda_N$ depends on $A$ and there may be less than $N$ negative eigenvalues. 

However in the latter case we can obviously replace $N$ by the lesser number $\bar{N}\Def \max (n\le N, \lambda_n\le 0)$ as $\E^*_N$ is decreasing function of $N$. In this case the first term  in (\ref{3-3}) would be just 
$\Tr^- (H_{A,W})$ and the second would be $0$. Then we apply theory of \cite{MQT11} immediately without extra complications.

Consider $A$ a minimizer (or its mollification) for potential $W=W^\TF$ and $\mu\le 0$. Then with an error $O(N^{\frac{2}{3}})$ 
\begin{equation}
\#\{ \lambda_k<\mu\} =\int (W -\mu )_+^{\frac{3}{2}}\, dx + O(N^{\frac{2}{3}}).
\label{3-4}
\end{equation}
One can prove (\ref{3-4}) easily using the regularity properties of $A$ established in \cite{MQT11} and the same rescaling as before. Note that the first term in (\ref{3-4}) differs from the same expression with $\mu =0$ (which is equal to $Z$) by $\asymp \mu (N^{4/3})^{1/2}\cdot N^{-1}=\mu N^{-1/3}$. Then as the left-hand expression equals $N$, and $N-Z=O(N^{\frac{2}{3}}$, we conclude that $|\lambda_N|=O(N)$.

Therefore modulo $O\bigl(N^{\frac{16}{9}}+\kappa a^{-3}N^2 \bigr)$ the sum of the first and the last term in (\ref{3-3}) is equal to
\begin{equation}
\frac{2}{15\pi^2} \int (W  -\lambda_N)_+^{\frac{5}{2}}\, dx + 
\sum _m 2 Z_m^2S(\kappa Z_m)
\label{3-5}
\end{equation}
and modulo $O(N^{-\frac{1}{3}}\lambda_N^2)=O(N^{\frac{5}{3}})$ one can rewrite the first term here as
\begin{equation}
\frac{2}{15\pi^2} \int W_+^{\frac{5}{2}}\, dx  -
\lambda_N   \frac{1}{3\pi^2} \int W_+^{\frac{3}{2}}\, dx
\label{3-6}
\end{equation}
and with the same error the second term here cancels term $\lambda_N N$ in (\ref{3-3}); then (\ref{3-3}) becomes 
\begin{equation}
\frac{2}{15\pi^2} 
\int W  _+^{\frac{5}{2}}\, dx +\sum _m 2Z_m^2S(\kappa Z_m)  -
\D(\rho, \rho)+ 
\D (\rho_\Psi-\rho , \rho_\Psi -\rho) 
\label{3-7}
\end{equation}
and as $W=W^\TF$, $\rho=\rho^\TF$ the first and the third term together are $\cE^\TF$, so we get again $\cE^\TF +\sum _m 2 Z_m^2S(\kappa Z_m) $.

Now we need to estimate properly the last term in (\ref{3-7}) i.e.
\begin{equation}
\frac{1}{2}\iint |x-y|^{-1}
\bigl(\rho_\Psi(x)-\rho^\TF (x)\bigr)\bigl(\rho_\Psi (y)-\rho^\TF(y)\bigr)\,dxdy. 
\label{3-8}
\end{equation}
Rescaling as before, and using (\ref{1-14}) we conclude that it does not exceed
\begin{equation}
N^{\frac{5}{3}}
\iint \varrho (x)^2\varrho(y)^2 \ell^{-1}(x)\ell^{-1}(y)|x-y|^{-1}\, dxdy
\label{3-9}
\end{equation}
where $\varrho$ is $\rho$ of \cite{MQT11} and we know that $\varrho=\ell^{-\frac{1}{2}}$ as $\ell\le 1$ and $\varrho =\ell^{-2}$ as $\ell\ge 1$.

Estimating integral by the (double) sum of integral as $\ell(x)\le 1$, $\ell(y)\le 1$ and $\ell(x)\ge 1$, $\ell(y)\ge 1$ we get (increasing $C$)
\begin{gather*}
C\int_{\{|y|\le |x|\le 1\}} |x-y|^{-1}|x|^{-2}|y|^{-2}\, dydx\asymp 1\\
\shortintertext{and}
C\int_{\{|y|\ge |x|\ge 1\}} |x-y|^{-1}|x|^{-3}|y|^{-3}\, dydx \asymp 1
\end{gather*}
respectively. 

This concludes the proof of the upper estimate in Theorem~\ref{thm-1-1} which is proven now.

\chapter{Miscellaneous}
\label{sect-4}

\begin{proof}{Proof of corollary \ref{cor-1-2}}
Optimization with respect to $\x_1,\ldots,\x_M$ implies 
\begin{equation}
\E^* < \sum_{1\le m\le M}\E^*_m
\label{4-1}
\end{equation}
where $\E^*=\E^*(\x_1,\ldots,\x_M;Z_1,\ldots,Z_m,N)$ and $\E^*_m=\E^*(Z_m, Z_m)$ are calculated for separate atoms. In virtue of theorem \ref{thm-1-1} and (\ref{1-9}) then
\begin{equation}
\cE^\TF-\sum_{1\le m\le M} \cE^\TF_m \le  Ca^{-3}N+ CN^{\frac{16}{9}};
\label{4-2}
\end{equation}
however due to strong non-binding theorem in Thomas-Fermi theory (see f.e. \cite{Sol}) the left-hand expression is $\asymp a^{-7}$ as 
$a\ge N^{-\frac{1}{3}}$ and therefore (\ref{4-2}) implies
\begin{equation}
a\ge \epsilon_1 N^{-\frac{16}{21}}
\label{4-3}
\end{equation}
and $a^{-3}N \le N^{\frac{16}{9}}$.

On the other hand, there is no binding with $a\le N^{-\frac{1}{3}}$ because remainder estimate is (better than) $CN^2$ and binding energy excess is $\asymp N^{\frac{7}{3}}$.
\end{proof}

\begin{remark}\label{rem-4-1}
Similar arguments work if we improve $N^{\frac{16}{9}}$ to $N^\nu$ with 
$\nu \ge \frac{7}{4}$ but without improving $a^{-3}N$ part of the remainder estimate we would not pass beyond $O(N^{\frac{7}{4}})$.
\end{remark}

There are several questions which after \cite{MQT11} could be answered in this framework by the standard arguments with certain error but we postpone it, hoping to improve remainder estimate $O(h^{-\frac{4}{3}})$ in  \cite{MQT11}:

\begin{problem}\label{problem-4-2}
\begin{enumerate}[label=(\roman*), leftmargin=*]
\item Investigate case $N\le Z-CZ^{\frac{2}{3}}$;
\item Estimate from above excess negative charge (how many extra electrons can and bind atom) ionization energy ($\E^*_{N-1}-\E^*_N$);
\item Estimate from above excess positive charge in the case of binding of several atoms i.e. estimate $Z-N$ as
\begin{equation}
\E^*(\x_1,\ldots,\x_M;Z_1,\ldots,Z_m,N) < \min _{\substack{N_1,\ldots,N_m: \\[2pt] N_1+\ldots +N_M=N}} \E^*_m (Z_m, N_m).
\label{4-4}
\end{equation}
\end{enumerate}
\end{problem}

\bibliographystyle{alpha}

\begin{thebibliography}{BrIvr}
\providecommand{\bysame}{\leavevmode\hbox to5em{\hrulefill}\thinspace}



\bibitem[BI]{bronstein:ivrii:IRO1}  M. Bronstein, M.,  V. Ivrii 
\href{http://weyl.math.toronto.edu:8888/victor2/preprints/IRO1.pdf}{Sharp Spectral Asymptotics for Operators with Irregular Coefficients. Pushing the Limits.}
\newblock \emph{Comm. Partial Differential Equations}, 28, 1\&2:99--123 (2003).


\bibitem[EFS1]{EFS1}  L. Erd{\H o}s, S. Fournais, J.P. Solovej: 
\emph{Stability and semiclassics in self-generated
fields.} 
\href{http://arxiv.org/abs/1105.0506}{arXiv:1105.0506}


\bibitem[EFS2]{EFS2}  L. Erd{\H o}s, S. Fournais,
J.P. Solovej: \emph{Second order semiclassics with self-generated
magnetic fields.} 
\href{http://arxiv.org/abs/1105.0512}{arXiv:http://arxiv.org/abs/1105.0512}


\bibitem[EFS3]{EFS3}  L. Erd{\H o}s, S. Fournais, J.P. Solovej: 
\emph{Scott correction for large atoms and molecules in a self-generated magnetic field} \href{http://arxiv.org/abs/1105.0521}{arXiv:1105.0521}


\bibitem[ES3]{ES3}  L. Erd{\H o}s, J.P. Solovej: 
\emph{Ground state energy of large atoms in a self-generated magnetic field} \href{http://arxiv.org/abs/0903.1816}{arXiv:0903.1816}
Commun. Math. Phys. \textbf{294}, No. 1, 229-249 (2009)


\bibitem[FS]{FS} C.~Fefferman and L.A.~Seco: \emph{On the energy of a
    large atom}, Bull.~AMS \textbf{23}, 2, 525--530 (1990).

\bibitem[FSW1]{FSW1} R. L. Frank, H. Siedentop, S. Warzel:
\emph{The ground state energy of heavy atoms: 
relativistic lowering of the leading energy correction.} 
Commun. Math. Phys.  \textbf{278}
   no. 2, 549--566 (2008)
   


\bibitem[FSW2]{FSW2} R. L. Frank, H. Siedentop, S. Warzel:
\emph{The energy of heavy atoms according to Brown and Ravenhall: the Scott correction.}
Doc. Math. \textbf{14}, 463--516 (2009).







\bibitem[FLL]{FLL}
J. Fr\"ohlich, E. H. Lieb, and M. Loss:
\emph{Stability of Coulomb systems with magnetic fields. 
I. The one-electron atom. }
Commun.\ Math.\ Phys.\ \textbf{104} 251--270 (1986) 

\bibitem[H]{H} W.~Hughes: \emph{An atomic energy bound that gives
    Scott's correction}, Adv.~Math. \textbf{79}, 213--270 (1990).


\bibitem[I1]{ivrii:MQT-Vishik}  V. Ivrii,
\href{http://weyl.math.toronto.edu:8888/victor2/preprints/MQT-Vishik.pdf}{Heavy molecules in the strong magnetic field.}
 Russian Journal of Math. Phys., 4(1):29--74 (1996).


   
\bibitem[I2]{ivrii:MQT1} V. Ivrii
\href{http://weyl.math.toronto.edu:8888/victor2/preprints/MQT1.pdf}
{Asymptotics of the ground state energy of heavy molecules in
the strong magnetic field. I}.  Russian Journal
of Mathematical Physics, 4(1):29--74 (1996).


\bibitem[I3]{ivrii:MQT2} V. Ivrii
\href{http://weyl.math.toronto.edu:8888/victor2/preprints/MQT2.pdf}
{Asymptotics of the ground state energy of heavy molecules in
the strong magnetic field. II}.  Russian Journal  of Mathematical Physics,
5(3):321--354 (1997).



\bibitem[I4]{ivrii:book} V. Ivrii
\href{http://books.google.com/books?id=AxChLigs2DAC&printsec=frontcover&dq=ivrii&ei=22rzR4iOFozsjgH99cWtDQ&sig=OkYod5fX163O8bIoK7kmg81atp0}{Microlocal Analysis and Precise Spectral Asymptotics}, Springer-Verlag, 1998.

\bibitem[I5]{ivrii:IRO2} V. Ivrii
\href{http://weyl.math.toronto.edu:8888/victor2/preprints/IRO2.pdf}{Sharp spectral asymptotics for operators with irregular coefficients.
Pushing the limits. II}. Comm. Part. Diff. Equats.,  28 (1\&2):125--156, (2003).




\bibitem[I6]{futurebook} V. Ivrii
 \emph{Microlocal Analysis and Sharp Spectral Asymptotics}, 
 in progress: available online at \newline
\href{http://www.math.toronto.edu/ivrii/futurebook.pdf}{http://www.math.toronto.edu/ivrii/futurebook.pdf}

 \bibitem[I7]{MQT10} V. Ivrii
 \emph{Local trace asymptotics  in the self-generated magnetic field}, 
 \newblock  \href{http://arxiv.org/abs/1108.4188}{arXiv:math/1108.4188 [math.AP]}:1–-24, 09/2011

 
 \bibitem[I8]{MQT11} V. Ivrii
 \emph{Global trace asymptotics in the self-generated magnetic field in the case of Coulomb-like singularities}, 
 \newblock  \href{http://arxiv.org/abs/1112.2487}{arXiv:math/1112.2487  [math.AP]}:1–-21, 12/2011.  

\bibitem[IS]{ivrii:ground} V. Ivrii  and , I.~M.~Sigal
\href{http://weyl.math.toronto.edu:8888/victor2/preprints/Scott.pdf}{Asymptotics of the ground state energies of large {C}oulomb systems.}
\newblock {Ann. of Math.}, 138:243--335 (1993).

\bibitem[L1]{Lieb-1} E. H. Lieb: \emph{Thomas-Fermi and related
theories of atoms and molecules}, Rev. Mod. Phys. \textbf{65}. No. 4, 603-641
(1981)



\bibitem[L2]{Lieb-2} E. H. Lieb: \emph{Variational principle for
many-fermion systems}, Phys. Rev. Lett. \textbf{46}, 457--459 (1981)
and \textbf{47} 69(E) (1981)




\bibitem[LLS]{LLS} E. H. Lieb, M. Loss and J. P. Solovej: \emph{ 
Stability of Matter in Magnetic Fields}, Phys. Rev. Lett. \textbf{75},
 985--989 (1995).


\bibitem[LO]{LO} E. H. Lieb and S. Oxford: \emph{Improved Lower Bound
    on the Indirect Coulomb Energy}, Int. J. Quant. Chem. \textbf{19},
  427--439, (1981)
  
\bibitem[S]{Sol} J.~P.~Solovej \emph{Asymptotic neutrality of diatomic molecules\/}, Commun. Math. Phys., \textbf{130}, 185-204, 1990.

 \end{thebibliography}



\end{document}